\def \NN {\mathbb N}
\def \CC {\mathbb C}

\def \RR {\mathbb R}
\def \ZZ {\mathbb Z}

\def \epsilon{\varepsilon}

\def \F  {{\mathcal F}}

\def \Q  {{\mathcal Q}}

\def \S  {{\mathcal S}}

\def \ep {\epsilon}
\def \ga {\gamma}

\def \si {\sigma}

\renewcommand{\S}{{\mathcal S}}

\documentclass[12pt,reqno]{amsart}

\usepackage{amsmath,amssymb}

\usepackage{url} 

\numberwithin{equation}{section}

\begin{document}


\title[]{Twists by Dirichlet characters and polynomial \\ Euler products of $L$-functions} 

\author[]{J.KACZOROWSKI \lowercase{and} A.PERELLI}
\maketitle

{\bf Abstract.} We prove that suitable properties of the twists by Dirichlet characters of an $L$-function of degree 2 imply that its Euler product is of polynomial type.



\vskip.5cm
\section{Introduction}

\smallskip
The properties of the twists by Dirichlet characters and the shape of the local Euler factors of $L$-functions are
two seemingly unrelated subjects. The main goal of the present paper is to show that in fact these two themes are closely related.  We focus our attention on  $L$-functions from the Selberg class $\S$ or the extended Selberg class $\S^{\sharp}$ as they provide a very convenient  framework for this study. We collect basic facts and notation related to $\S$ and  $\S^{\sharp}$ in Section 2. Here we recall that the twist of a function
\begin{equation}
\label{eq:1}
F(s) = \sum_{n=1}^{\infty}\frac{a(n)}{n^s}
\end{equation}
from $\S^{\sharp}$ by a Dirichlet character $\chi$ (mod $q$) is defined as
\[ 
F^{\chi}(s) = \sum_{n=1}^{\infty}\frac{a(n)\chi(n)}{n^s}.
\] 
This operation is fundamental in the theory of automorphic $L$-functions. If $F$ belongs to the Selberg class, it is expected that $F^\chi$ also belongs to the same class, at least when the conductors of $F$ and of the primitive character $\chi$ are coprime. This is true in many special cases, particularly for automorphic $L$-functions, but the general problem is wide open. 

\smallskip
The second theme is the study of the admissible shape of local Euler factors of $F$  in (\ref{eq:1}).  Recall that the local factor of $F$ at a prime $p$ is defined as
\begin{equation}
\label{eq:Fp}
 F_p(s) = \sum_{k=0}^{\infty} \frac{a(p^k)}{p^{ks}}.
 \end{equation}
The Euler product axiom in the definition of the Selberg class tells us that
\begin{equation}
\label{logFp}
F(s)=\prod_p F_p(s) \quad \text{and} \quad \log F_p(s) = \sum_{k=1}^{\infty}\frac{b(p^k)}{p^{ks}},
 \end{equation}
where
\begin{equation}
\label{bpk}
 b(p^k)\ll p^{\vartheta k}
 \end{equation}
for a certain $\vartheta <1\slash 2$. It is expected that every $F\in S$ has a polynomial Euler product, i.e. of type
\[ 
F_p(s) = \prod_{j=1}^{\partial_p} \left( 1-\frac{\alpha_{j,p}}{p^s}\right)^{-1}
\]
for all primes $p$. Again, this is true in many special cases, in particular for automorphic $L$-functions, but in general the problem is wide open.

\smallskip
In this paper we show, in the case of functions of degree 2, that suitable properties of the twists by Dirichlet characters imply that the local factors are of polynomial type. In order to state our results in a synthetic way we adopt the following terminology.

\smallskip
Given a prime $p$, we say that $F\in\S^\sharp$ {\it splits at} $p$ if for $\si>1$
\[
F(s) = F_p(s) \sum_{p\nmid n} \frac{a(n)}{n^s}
\]
and $F_p(s)$ satisfies (\ref{logFp}) and (\ref{bpk}). Note that if $F$ splits at $p$ then $a(1)=1$ and $a(n) = a(p^\ell) a(k)$ whenever $n=p^\ell k$ with $p\nmid k$. In particular, an $L$-function from the Selberg class splits at all primes $p$. Let $F\in\S^\sharp$, $q_F$ be its conductor, assume that $q_F\in\NN$ and let $p$ be a prime not dividing $q_F$.  We denote by $m_{q_F}(p)$ the order of $p$ (mod $q_F$), i.e. the least positive integer $m$ such that $p^m\equiv 1$ (mod $q_F$). We say that $F\in\S^\sharp$ is {\it weakly twist-regular} at $p$ if for every primitive Dirichlet character $\chi$ (mod $p^f$) with $1\leq f\leq m_{q_F}(p)$, the twist $F^\chi$ belongs to $\S^\sharp$ and has the same degree as $F$. Moreover, we say that $F\in\S^\sharp$ is {\it twist-regular} at $p$ if it is weakly twist-regular at $p$ and for $f=1$ the conductor $q_\chi$ of $F^\chi$ satisfies $q_\chi=q_Fp^{d_F}$, $d_F$ being the degree of $F$.

\medskip
{\bf Theorem 1.} {\sl Let $F\in\S^\sharp$ be of degree $2$ and conductor $q_F\in\NN$. Then there exists a constant $B_F>0$ such that if $F$ splits and is weakly twist-regular at a prime $p>B_F$, then its local factor $F_p(s)$ is a rational function of $p^{-s}$ and its numerator has degree $\leq m_{q_F}(p)-1$. In particular, if $p\equiv 1$  {\rm (mod $q_F$)} then $F_p(s)^{-1}$ is a polynomial in $p^{-s}$. Moreover, if $F\in\S$ then the condition ``$p>B_F$'' can be replaced by ``$p \nmid q_F$''.}

\medskip
{\bf Theorem 2.} {\sl Let $F\in\S^\sharp$ be of degree $2$ and conductor $q_F\in\NN$. Then there exists a constant $B_F>0$ such that if $F$ splits and is weakly twist-regular at two distinct primes $p,q>B_F$, $p\equiv q$ {\rm (mod $q_F$)}, then $F_p(s)^{-1}$ and $F_q(s)^{-1}$ are polynomials in $p^{-s}$ and $q^{-s}$, respectively. Moreover, if $F\in\S$ then the condition ``$p,q>B_F$'' can be replaced by ``$p,q \nmid q_F$''.}

\medskip
An interesting feature of the L-functions, depending on the so-called multiplicity one property, is that the local factors at the primes not dividing the conductor determine all the other, i.e. those at the ramified primes. The next theorem makes this phenomenon more precise for prime conductors $q_F$. Indeed, we show that the local factor at $q_F$ is of polynomial type as well. We expect that a similar result holds for general $q_F$, but it seems that proof of this fact would require a non-trivial extension of the methods in the present paper. We hope to address this problem in a future paper.

\medskip
{\bf Theorem 3.} {\sl Let $F\in\S$ be of degree $2$ and its conductor $q_F$ be a prime. If $F$ is weakly twist-regular at all primes $\neq q_F$, then $F$ has a polynomial Euler product.}

\medskip
Our final result gives information on the degree of the polynomials $F_p(s)^{-1}$.

\medskip
{\bf Theorem 4.} {\sl Let $F\in\S$ be of degree $2$ and conductor $q_F\in\NN$. If $F$ is twist-regular at all primes $p\nmid q_F$, then there exists a constant $B_F>0$ such that for all primes $p\geq B_F$ we have
\[
F_p(s) = \left(1-\frac{\alpha_p}{p^s}\right)^{-1} \left(1-\frac{\beta_p}{p^s}\right)^{-1}
\]
with certain $|\alpha_p|,|\beta_p|\leq 1$. Moreover, for the primes $p<B_F$, $p\nmid q_F$, we have that $F_p(s)^{-1}$ is a polynomial in $p^{-s}$.}

\medskip
We finally remark that the results in this paper are a step toward a very general form of Weil's converse theorem for $L$-functions of degree 2. We shall address this problem in forthcoming papers.
 
\medskip
{\bf Acknowledgements}.  This research was partially supported by the Istituto Nazionale di Alta Matematica, by the MIUR grant PRIN-2017 {\sl ``Geometric, algebraic and analytic methods in arithmetic''} and by grant 2021/41/BST1/00241 {\sl ``Analytic methods in number theory''}  from the National Science Centre, Poland.

\section{Notation}

\smallskip
Throughout the paper we write $s=\si+it$, $e(x) = e^{2\pi ix}$ and $\overline{f}(s)$ for $\overline{f(\overline{s})}$. The extended Selberg class $\S^\sharp$ consists of non identically vanishing Dirichlet series \eqref{eq:1} absolutely convergent for $\si>1$, such that $(s-1)^mF(s)$ is entire of finite order for some integer $m\geq0$, and satisfying a functional equation of type
\[
F(s) \gamma(s) = \omega \overline{\gamma}(1-s) \overline{F}(1-s),
\]
where $|\omega|=1$ and the $\gamma$-factor
\[
\gamma(s) = Q^s\prod_{j=1}^r\Gamma(\lambda_js+\mu_j) 
\]
has $Q>0$, $r\geq0$, $\lambda_j>0$ and $\Re(\mu_j)\geq0$. Note that the conjugate function $\overline{F}$ has conjugated coefficients $\overline{a(n)}$. The Selberg class $\S$ is the subclass of $\S^\sharp$ of the functions with an Euler product as in \eqref{eq:Fp},\eqref{logFp} and \eqref{bpk}, and satisfying the Ramanujan conjecture $a(n) \ll n^\ep$. Note that the local factors $F_p(s)$ in \eqref{eq:Fp} satisfy
\begin{equation}
\label{effepi}
F_p(s)\neq0 \quad \text{for $\si>\vartheta$.}
\end{equation}
We refer to our survey papers \cite{Kac/2006},\cite{Ka-Pe/1999b},\cite{Per/2005},\cite{Per/2004},\cite{Per/2010},\cite{Per/2017} for further definitions, examples and the basic theory of the Selberg class. 

\smallskip
Degree $d_F$, conductor $q_F$, root number $\omega_F$ and $\xi$-invariant $\xi_F$ of $F\in\S^\sharp$ are defined by
\[
\begin{split}
d_F =2\sum_{j=1}^r\lambda_j, \qquad q_F= (2\pi)^dQ^2\prod_{j=1}^r\lambda_j^{2\lambda_j}, \\
\omega_F=\omega \prod_{j=1}^r \lambda_j^{-2i\Im(\mu_j)},  \qquad \xi_F = 2\sum_{j=1}^r(\mu_j-1/2)= \eta_F+ id_F\theta_F
\end{split}
\]
with $\eta_F,\theta_F\in\RR$. We also write
\[
\omega_F^* = \omega_F e^{-i\frac{\pi}{2}(\eta_F+1)}\big(\frac{q_F}{(2\pi)^{2}}\big)^{i\frac{\theta_F}{2}} \quad \text{and} \quad \tau_F=\max_{j=1,\dots,r}\big|\frac{\Im(\mu_j)}{\lambda_j}\big|,
\]
while $m_F$ denotes the order of the pole of $F$ at $s=1$.

\smallskip
Finally, the linear twist of $F\in\S^\sharp$ is defined as
\[
F(s,\alpha) = \sum_{n=1}^\infty \frac{a(n)}{n^s} e(-n\alpha)
\]
with $\alpha\in\RR$.

\section{Lemmas}

\smallskip
Given a Dirichlet character $\chi$ we denote, as usual, by $\chi^*$ the primitive character inducing $\chi$ and by $\chi_0$ the principal character. Also, we recall that $\mu(n)$ and $\varphi(n)$ denote the M\"obius and Euler functions, respectively.

\medskip
{\bf Lemma 1.} {\it Let $p$ be a prime number. For every integer $r\geq1$ there exist coefficients $c(\chi,p^r)$, where $\chi$ runs over the Dirichlet characters} (mod $p^r$), {\it such that for $(n,p^r)=1$ we have
\[
e(-n/p^r) = \sum_{\chi \, {\rm (mod} \, p^r {\rm )}} c(\chi,p^r) \chi(n).
\]
Moreover
\[
c(\chi_0,p^r) = 
\begin{cases}
\frac{1}{1-p} & \text{if} \ r=1 \\
0 & \text{if} \ r>1.
\end{cases}
\]}

\medskip
{\it Proof.} The existence of the coefficients $c(\chi,p^r)$ follows by elementary harmonic analysis on the group $\ZZ_{p^r}^*$ of the reduced residues (mod $p^r$). Moreover, by orthogonality we have that
\[
c(\chi_0,p^r) = \frac{1}{\varphi(p^r)} \sum_{n\in\ZZ_{p^r}^*} e(-n/p^r) = \frac{\mu(p^r)}{\varphi(p^r)},
\]
and the lemma follows. \qed

\medskip
Given a function $F\in\S^\sharp$, a prime $p$ and a positive integer $m$ we define the polynomial $W_{m,p}(X)$ as
\[
W_{m,p}(X) = \sum_{\ell=0}^{m-2} a(p^\ell) X^\ell + \frac{p}{p-1} a(p^{m-1}) X^{m-1}.
\]
Note that if $m=1$ the sum is empty and hence equals 0. Note also that $W_{m,p}(X)$ is not identically vanishing since $a(1)=1$, as $F$ splits at $p$.

\medskip
{\bf Lemma 2.} {\it Let $p$ be a prime number, $m$ be a positive integer and suppose that $F\in\S^\sharp$ splits at $p$. Then for $\si>1$ we have}
\[
F(s,1/p^m) = \sum_{\ell=0}^{m-1} \frac{a(p^\ell)}{p^{\ell s}} \sum_{\substack{\chi \, {\rm (mod} \, p^{m-\ell} {\rm)} \\ \chi\neq \chi_0}} c(\chi,p^{m-\ell}) F^{\chi^*}(s) + \left(1- W_{m,p}(p^{-s}) F_p(s)^{-1}\right) F(s).
\]

\medskip
{\it Proof.} Writing $n=p^\ell k$ with $p\nmid k$, since $F$ splits at $p$ for $\si>1$ we have
\[
\begin{split}
 \sum_{n=1}^{\infty} \frac{a(n)}{n^s} e(-n/p^m) &=  \sum_{\ell=1}^{m-1} \sum_{p^{\ell}\| n} \frac{a(n)}{n^s} e(-n/p^m)
 +\sum_{p^m|n} \frac{a(n)}{n^s} +\sum_{p\nmid n} \frac{a(n)}{n^s} e(-n/p^m) \\
 &= \sum_{\ell=0}^{m-1}\frac{a(p^{\ell})}{p^{\ell s}} \sum_{p\nmid k} \frac{a(k)}{k^s} e(-k/p^{m-\ell})
 + F(s) - \sum_{\ell=0}^{m-1}\frac{a(p^{\ell})}{p^{\ell s}}\sum_{p\nmid k} \frac{a(k)}{k^s}.
 \end{split}
\]
Now we apply Lemma 1 to obtain that
\[
\begin{split}
\sum_{p\nmid k} \frac{a(k)}{k^s} e(-k/p^{m-\ell}) &= \sum_{\chi \, {\rm (mod} \, p^{m-\ell} {\rm )}} c(\chi,p^{m-\ell}) F^\chi(s) \\
&=  \sum_{\substack{\chi \, {\rm (mod} \, p^{m-\ell} {\rm )} \\ \chi\neq\chi_0}} c(\chi,p^{m-\ell}) F^{\chi^*}(s) + c(\chi_0,p^{m-\ell}) F^{\chi_0}(s).
\end{split}
\]
Moreover, since $F^{\chi_0}(s) = \sum_{p\nmid n} a(n) n^{-s}= F_p(s)^{-1}F(s)$, from Lemma 1 we also have that
\[
c(\chi_0,p^{m-\ell}) F^{\chi_0}(s) =
\begin{cases}
\frac{1}{1-p} F_p(s)^{-1} F(s) & \text{if} \ \ell=m-1 \\
0 & \text{if} \ \ell<m-1,
\end{cases}
\]
and the lemma follows by a simple computation. \qed

\medskip
{\bf Lemma 3.} {\it Let $F\in\S^\sharp$ with $d=2$, and let $a<b$ be fixed. Then
\[
F(s) \ll \left(\frac{q_F}{(2\pi e)^2}\right)^{|\si|} |\si|^{2|\si|+1}
\]
uniformly for $a\leq t\leq b$ and $\si\leq-1$. Moreover, if $[a,b]\cap[-\tau_F,\tau_F]=\emptyset$ we also have
\[
F(s) \gg \left(\frac{q_F}{(2\pi e)^2}\right)^{|\si|} |\si|^{2|\si|+1}
\]
uniformly for $a\leq t\leq b$ and $\si\leq-1$.}

\medskip
{\it Proof.} This is a slightly refined version of Lemma 2.1 in \cite{Ka-Pe/2015}. We follow the proof of that lemma, using the notation there and recalling that here we have $d=2$. Then we observe that Stirling's formula actually gives the more precise expression
\[
\log|G(s)| = 2\si\log\si + (\log\beta-2)\si +\log\si +O(1),
\]
and the lemma follows. \qed

\medskip
{\bf Lemma 4.} {\it Let $F\in\S^\sharp$ with $d_F=2$, and let $\alpha>0$. Then for every integer $K>0$ there exist polynomials $Q_0(s),...,Q_K(s)$, with $Q_0(s)\equiv 1$, such that
\begin{equation}
\label{eq:Falpha} 
F(s,\alpha) = -i\omega_F^* (\sqrt{q_F}\alpha)^{2s-1+i\theta_F} \sum_{\nu=0}^K \big(\frac{iq_F\alpha}{2\pi}\big)^\nu Q_\nu(s) \overline{F}\big(s+\nu+2i\theta_F,-\frac{1}{q_F\alpha}\big) + H_K(s,\alpha).
\end{equation}
Here $H_K(s,\alpha)$ is holomorphic for $-K+\frac12<\sigma<2$ and $|s|<2K$, and satisfies
\begin{equation}
\label{HK}
H_K(s,\alpha) \ll (AK)^K
\end{equation}
with a certain constant $A=A(F,\alpha)>0$. Moreover, $\deg Q_\nu=2\nu$ and}
\begin{equation}
\label{eq:Q1}
Q_\nu(s) \ll \frac{(A(|s|+1))^{2\nu}}{\nu!} \hskip1.5cm \text{for} \ 1\leq \nu\leq \min(|s|,K)
\end{equation}
\begin{equation}
\label{eq:Q2}
Q_\nu(s) \ll (AK)^K \hskip1.5cm \text{for} \  |s|\leq 2K, \, \nu\leq K.
\end{equation}

\medskip
{\it Proof.} This is Theorem 2.1 of \cite{Ka-Pe/2015}. Note that the different value of the shift in \eqref{eq:Falpha} with respect to Theorem 2.1, namely $\nu+2i\theta_F$ in place of $\nu+i\theta_F$, is due to the slightly different definition of $\theta_F$ used in this paper. \qed

\medskip
{\bf Lemma 5.} {\it Let $a<b$ be fixed and let ${\Q}_{\nu}(s),{\F}(s)$ be real-valued continuous functions defined in the horizontal strip $a\leq t\leq b$ satisfying the following conditions:
\begin{equation}
\label{F1}
{\F}(s)\ll 1\ \hskip1.5cm \text{for}\ \si\geq1,
\end{equation}
\begin{equation}
\label{F2}
{\F}(s)\ll B^{|\sigma|}|\sigma|^{2|\sigma|+1} \hskip1.5cm \text{for}\ \ \sigma\leq 1,
\end{equation}
\begin{equation}
\label{mQ1}
 {\Q}_{\nu}(s) \ll \frac{(C(|s|+1))^{2\nu}}{\nu!}  \hskip1.5cm \text{for}\ \ 0\leq \nu\leq \min(|s|,|\sigma|+2),
 \end{equation}
\begin{equation}
\label{mQ2}
{\Q}_{\nu}(s) \ll (C(|\sigma|+2))^{|\sigma|+2} \hskip1.5cm \text{for}\ \ \nu\leq |\sigma|+2\ ,\ |s|\leq 2[|\sigma|]+2,
\end{equation}
 where $B,C>0$.  Then for $\sigma\leq -1$ and $a\leq t\leq b$ we have
 \begin{equation}
 \label{Phi}
\sum_{0\leq\nu\leq|\sigma|+2} {\mathcal Q}_{\nu}(s){\mathcal F}(s+\nu) \ll B^{|\sigma|}|\sigma|^{2|\sigma|+1}
 \end{equation}
with the implied constant depending on $\F,a,b,B,C$ and implied constants in \eqref{F1}-\eqref{mQ2}.}
 
 \medskip
 {\it Proof.} In the proof we use the synthetic expression ``suitably bounded'' to mean ``bounded by a constant depending at most on $\F,a,b,B,C$ and implied constants in \eqref{F1}-\eqref{mQ2}''. We first note that we may assume without loss of generality that
\begin{equation}
\label{sigma}
|\sigma|\geq \max(|a|,|b|) +1.
\end{equation}
Indeed, otherwise both $s$ and $s+\nu$ with $0\leq\nu\leq|\sigma|+2$ stay in a compact domain and hence the functions ${\Q}_{\nu}(s),{\F}(s+\nu)$ are suitably bounded. Thus $\Phi(s)$ is suitably bounded as well and the assertion follows in this case.

\smallskip
Assuming \eqref{sigma}, we first  consider the terms in \eqref{Phi} with
\begin{equation}
\label{nu1}
|\sigma|-\max(|a|,|b|)-1\leq \nu\leq |\sigma|+2.
\end{equation}
For such $\nu$ we have that $-\max(|a|,|b|)-1\leq -|\sigma| +\nu\leq 2$; thus, as before, $s+\nu$ stays in a compact domain and hence ${\F}(s+\nu)$ is suitably bounded. Moreover, recalling (\ref{sigma}) we certainly have
\[
|s|\leq |\sigma|+|t|\leq 2[|\sigma|]+2,
\]
so in view of this and of \eqref{nu1} we can apply \eqref{mQ2} to estimate ${\Q}_{\nu}(s)$. Hence the terms in $\Phi(s)$ corresponding to $\nu$ in the range \eqref{nu1} contribute at most
\[
\ll (C(|\sigma|+2))^{|\sigma|+2} \ll B^{|\sigma|} |\sigma|^{2|\sigma|+1},
\]
and our assertion holds in this case as well.

\smallskip
Finally suppose that
\begin{equation}
\label{nu2}
0\leq \nu \leq |\sigma|-\max(|a|,|b|)-1.
\end{equation}
Recalling that $\si$ is negative, for such $\nu$ we have $\sigma+\nu\leq -1$ and hence we can apply \eqref{F2} to get
\[
\F(s+\nu) \ll |\sigma| B^{|\sigma|-\nu} (|\sigma|-\nu)^{2(|\sigma|-\nu)}.
\]
Moreover, for such $\nu$ we also have $\nu\leq \min(|s|,|\si|+2)$ so we can apply \eqref{mQ1} to obtain, thanks to \eqref{sigma}, that
\[
{\Q}_{\nu}(s)\ll \frac{(C(|s|+1))^{2\nu}}{\nu!} \leq \frac{(2C|\sigma|)^{2\nu}}{\nu!}.
\]
Hence the terms corresponding to $\nu$ in the range \eqref{nu2} contribute at most
\[
\begin{split}
&\ll |\sigma|\sum_{0\leq\nu\leq|\sigma|} \frac{(2C|\sigma|)^{2\nu}}{\nu!} B^{|\sigma|-\nu} (|\sigma|-\nu)^{2(|\sigma|-\nu)} \\
 &\ll |\sigma|B^{|\sigma|}
 \max_{0\leq \nu\leq |\sigma|} \left(|\sigma|^{2\nu} (|\sigma|-\nu)^{2(|\sigma|-\nu)}\right) \sum_{\nu=0}^{\infty}\frac{1}{\nu!}\left(\frac{4C^2}{B}\right)^{\nu} \\
 &\ll B^{|\sigma|} |\sigma|^{2|\sigma|+1},
\end{split}
\]
and the lemma follows. \qed

\medskip
{\bf Lemma 6.} {\it  Let $F\in S^{\sharp}$ with $d=2$. Then the linear twist $F(s, 1/q_F)$ has meromorphic continuation to $\CC$ with poles at most at the points
 \[ 
 1-\nu - 2i\theta_F \hskip1.5cm \nu\in\ZZ, \nu\geq0
 \]
and order $\leq m_F$. Moreover, for every $a<b$ such that
\begin{equation}
\label{eq:ab}
 (a+\theta_F)(b+\theta_F)>0
 \end{equation}
 we have, uniformly for $\si\leq -1$ and $a\leq t \leq b$, that}
 \begin{equation}
 \label{F1/q}
 F(s,1/q_F) \ll_{F,a,b}  \left(\frac{q_F}{2\pi e}\right)^{2|\sigma|} |\sigma|^{2|\sigma|+1}.
 \end{equation}

 \medskip
 {\it Proof.} We start with Lemma 4 with the choice $\alpha=1/q_F$. Since
\begin{equation}
\label{F1/qa}
 \overline{F}\big(s+\nu+2i\theta_F, -\frac{1}{q_F\alpha}\big) =\overline{F}( s+\nu+2i\theta_F),
\end{equation} 
the terms on the right hand side of \eqref{eq:Falpha} are holomorphic for $s\neq 1-\nu-2i\theta_F$, with $\nu\geq 0$ and $s$ in the range stated after \eqref{eq:Falpha}. Moreover, the potential poles are induced by $F(s)$ and hence of order $\leq m_F$. The first assertion follows since $K$ is arbitrarily large.

\smallskip
To prove the boud \eqref{F1/q} we use Lemma 4 with $K=[|\sigma|]+2$. Recalling  \eqref{F1/qa}, from \eqref{eq:Falpha} and \eqref{HK} we have that
\begin{equation}
\label{F1/qa2}
F(s,1/q_F) \ll q_F^{|\sigma|}\sum_{0\leq \nu\leq|\sigma|+2} (2\pi)^{-\nu} |Q_{\nu}(s)| |\overline{F}(s+\nu+2i\theta_F)| + O\big((A(|\sigma|+2))^{|\sigma|+2}\big).
\end{equation}
Next we apply Lemma 5 with 
\[ 
\Q_{\nu}(s)= (2\pi)^{-2\nu}|Q_{\nu}(s)| \quad \text{and} \quad \F(s)=  |\overline{F}(s+i\theta)|.
\]
From \eqref{eq:Q1} and \eqref{eq:Q2} we see that \eqref{mQ1} and \eqref{mQ2} hold with $C=A$, and from Lemma 3 we see that \eqref{F2} holds with $B=q_F/(2\pi e)^2$. Finally, \eqref{F1} is satisfied as well, thanks to \eqref{eq:ab} and the description of the possible singularities of $F(s,1/q_F)$. Now it is clear that \eqref{F1/q} follows from \eqref{Phi} and \eqref{F1/qa2}, and the proof is complete. \qed

\medskip
For $F$ in the Selberg class $\S$ we denote by $N_F(\si,T)$ the number of non-trivial zeros $\beta+i\gamma$ of $F$ in the rectangle $\beta> \si$, $|\ga|\leq T$.

\medskip
{\bf Lemma 7.} {\it Let $F\in\S$ with $d=2$. Then for every $\ep>0$ and any fixed $\si>1/2$ we have}
\[
N_F(\si,T) \ll T^{3/2-\si+\ep}.
\]

\medskip
{\it Proof.} See p.474-475 of \cite{Ka-Pe/2015}. \qed

\section{Proof of Theorem 1}

\smallskip
Let $p$ be as in Theorem 1. In particular we may assume that $p\nmid q_F$, and let $m=m_{q_F}(p)$. According to Lemma 2 we have
\begin{equation}
\label{T1.1}
W_{m,p}(p^{-s}) F_p(s)^{-1} F(s) = \sum_{\ell=0}^{m-1} \frac{a(p^{\ell})}{p^{\ell s}} \sum_{\substack{\chi \, {\rm (mod} \, p^{m-\ell} {\rm)} \\ \chi\neq \chi_0}} c(\chi, p^{m-\ell}) F^{\chi^*}(s)
+ F(s)- F(s,1/p^m),
\end{equation}
where all terms on the right-hand side are meromorphic on $\CC$ except possibly the last one. Since $p^m\equiv 1$ (mod $q_F$) we have $F(s,p^m/q_F) = F(s, 1/q_F)$, hence from Lemma 4 with $\alpha=1/p^m$ we obtain
\begin{equation}
\label{T1.2}
\begin{split}
F(s, 1/p^m) =&-i\omega_F^* \left(\frac{\sqrt{q_F}}{p^m}\right)^{2s-1+i\theta_F}
\sum_{\nu=0}^K\left(\frac{iq_F}{2\pi p^m}\right)^{\nu} Q_{\nu}(s) \overline{F} (s+\nu+2i\theta_F, -1/q_F) \\
 &+ H_K(s,1/p^m).
 \end{split}
\end{equation}
By Lemma 6 this gives meromorphic continuation of $F(s,1/p^m)$ to $\CC$, and hence the same is true for the term on the left hand side of \eqref{T1.1}. Therefore, the function
\begin{equation}
\label{WFp}
W(s) := W_{m,p}(p^{-s}) F_p(s)^{-1}
\end{equation}
 is meromorphic on $\CC$ and, thanks to \eqref{effepi}, its possible poles lie on the half-plane $\sigma\leq\vartheta$. Recalling \eqref{T1.1},\eqref{T1.2} and the structure of the singularities of $F(s,1/q_F)$ described in Lemma 6, the possible poles of $W(s)$ outside a certain horizontal strip of bounded height are induced by the zeros of $F(s)$. But $F_p(s)$ is $(2\pi i/\log p)$-periodic, hence one such pole generates $\gg T$ poles in the strip $\si\leq \vartheta$, $|t|\leq T$. If $F\in\S$ this contradicts the density estimate in Lemma 7, thus $W(s)$ is entire. If $F\in\S^\sharp$ we may apply an idea from Gierszewski \cite{Gie/2020} and choose $B_F>0$ so large that the whole horizontal strip $2\pi/\log B_F\leq t\leq 4\pi/\log B_F$ is free from non-trivial zeros of $F$, thus deducing that $W(s)$ is entire if $p>B_F$.
 
\smallskip
To conclude the proof it suffices to show that $W(s)$ is a polynomial in $p^{-s}$, since $W_{m,p}(p^{-s})$ is a polynomial in $p^{-s}$ of degree $\leq m-1$. To this end, since $W(s)$ is $(2\pi i/\log p)$-periodic, we start by estimating of $W(s)F(s)$ for $a\leq t\leq b$, where $a>\max(|\theta_F|,\tau_F)$ and $b=a+2\pi/\log p$. From \eqref{T1.1} we obtain
\begin{equation}
\label{T1.3}
W(s) F(s) \ll C_1^{|\si|} \max_{0\leq\ell\leq m-1} \max_{\substack{\chi \, {\rm (mod} \, p^{m-\ell} {\rm)} \\ \chi\neq \chi_0}} |F^{\chi^*}(s)| + |F(s)|+ |F(s, 1/p^m)|
\end{equation}
with a certain constant $C_1>0$. Since $F$ is weakly twist-regular at $p$ we have that $F^{\chi^*}\in\S^\sharp$ and has degree 2, hence by Lemma 3 the sum of the first two terms in \eqref{T1.3} is 
\[
\ll C_2^{|\sigma|}|\sigma|^{2|\sigma|+1}
\] 
with $C_2>0$.
Moreover, from \eqref{T1.2} and Lemma 4 with $K=[|\sigma|]+2$ we see that the last term is
\[ 
\ll \left(\frac{p^m}{\sqrt{q_F}}\right)^{2|\sigma|} \sum_{\nu=0}^K \left(\frac{q_F}{2\pi p^m}\right)^{\nu} |Q_{\nu}(s)| 
|F(\overline{s}+\nu-2i\theta_F,1/q_F)|  +O \left( (A|\sigma|)^{|\sigma|}\right).
\]
We estimate this sum using Lemma 5 with
\[
\Q_{\nu}(s)=  \left(\frac{q_F}{2\pi p^m}\right)^{\nu} |Q_{\nu}(s)|  \quad \text{and} \quad \F(s) = |F(\overline{s}-2i\theta_F, 1/q_F)|.
\]
As in the proof of Lemma 6, but using Lemma 6 itself in place of Lemma 3,  we easily check that the assumptions in Lemma 5 are satisfied, hence concluding that
\[ 
F(s, 1\slash p^m)\ \ll C_3^{|\sigma|} |\sigma|^{2|\sigma|+1}
\] 
with some $C_3>0$.
Thus
\begin{equation}
\label{WFF}
W(s) F(s) \ll C_4^{|\sigma|} |\sigma|^{2|\sigma|+1}
\end{equation}
with $C_4=\max(C_2,C_3)$. 
On the other hand, we estimate $F(s)$ from below using Lemma 3. Recalling \eqref{WFF} and that $a>\tau_F$ we obtain
\[
W(s)\ll C_5^{|\sigma|}
\]
with a certain $C_5=C_5(F,p)>0$, uniformly for $-\infty< \sigma <\infty$ and $a\leq t\leq b$. As in the proof of Theorem 1.1 of \cite{Ka-Pe/2015} (see p.448), this bound implies that $W(s)$ is a polynomial in $p^{-s}$, as required.

\smallskip
Finally, we write
\[ 
F_p(s) = \frac{{N}_p(p^{-s})}{{D}_p(p^{-s})}
\]
for certain coprime normalized polynomials ${N}_p, {D}_p\in {\mathbb C}[X]$. Since $W(s)$ is entire, from \eqref{WFp} we see that ${N}_p$ divides $W_{m,p}$, thus has degree $\leq m-1$. In particular, $N_p\equiv 1$ if $p\equiv1$ (mod $q_F$), and Theorem 1 follows. \qed

\section{Proof of Theorem 2}

\smallskip
From Theorem 1 we know that $F_p(s)^{-1}$ is a rational function of $p^{-s}$; suppose that this rational function  is not a polynomial. Hence the singularities of $F_p(s)^{-1}$ lie on a finite number of vertical lines and, in view of \eqref{effepi}, there exists $\sigma_p\leq \vartheta <1\slash 2$ such that  $F_p(s)^{-1}$ is holomorphic for $\sigma>\sigma_p$ and has infinitely many poles on the line $\sigma=\sigma_p$. More precisely, the poles of $F_p(s)^{-1}$ lie on finitely many arithmetic progressions with difference $2\pi i\slash \log p$ and, thanks to Lemma 7, every such arithmetic progression contains infinitely many poles of $F_p(s)^{-1}F(s)$. Moreover, applying Lemma 2 with $m=1$ we obtain that
\begin{equation}
\label{t2-1}
\frac{p}{p-1}F_p(s)^{-1} F(s) = \sum_{\substack{\chi \, {\rm (mod} \, p {\rm)} \\ \chi\neq \chi_0}} c(\chi, p)F^{\chi}(s) +F(s) - F(s,1/p),
\end{equation}
hence $F(s,1/p)$ is meromorphic over $\CC$, has the same singularities of $F_p(s)^{-1}F(s)$ on the line $\si=\si_p$ and at most a pole at $s=1$ for $\si>\si_p$. It is also clear that analogous statements are true with the prime $q$ in place of $p$, assuming that $F_q(s)$ is not of polynomial type, and without loss of generality  we may assume that $\sigma_q\leq \sigma_p$. If $F_q(s)$ is of polynomial type we set $\sigma_q=-\infty$.

\smallskip
We denote by $h(s)$ a generic function which is meromorphic with finitely many poles for $\sigma>\sigma_p-1$. From Lemma 4 with $\alpha= p/q_F$ and arbitrarily large $K$, and observing that the terms with $\nu$ sufficiently large are holomorphic for $\sigma>\sigma_p-1$, we have
\[
F(s,p/q_F) = -i\omega_F^* (p/\sqrt{q_F})^{2s-1+i\theta_F} \sum_{\nu=0}^{K_p} \left(\frac{ip}{2\pi}\right)^\nu Q_\nu(s) \overline{F}(s+\nu+2i\theta_F,-1/p) + h(s)
\]
with a suitable integer $K_p$. Hence, recalling that $Q_0(s)\equiv1$, from the above properties of $F(s,1/p)$ we deduce that
\begin{equation}
\label{5-2}
F(s,p/q_F) = -i\omega_F^* (p/\sqrt{q_F})^{2s-1+i\theta_F} \overline{F}(s+2i\theta_F,-1/p) + h(s).
\end{equation}
Again, an analogous formula holds with the prime $q$ in place of $p$. Moreover, by a double conjugation from \eqref{5-2} we obtain that
\begin{equation}
\label{eq:t2-1a}
F(s,1/p) =-i\omega_F^* (\sqrt{q_F}/p)^{2s-1+i\theta_F} \overline{F(\overline{s}-2i\theta_F,p/q_F)} +  h(s),
\end{equation}
and similarly with $q$ in place of $p$.

\smallskip
But $p\equiv q$ (mod $q_F$), therefore
\[
F(\overline{s}-2i\theta,p/q_F)=F(\overline{s}-2i\theta_F,q/q_F)
\]
and hence from \eqref{eq:t2-1a} we obtain that
\[
F(s,1/p) = (q/p)^{2s-1+i\theta_F} F(s,1/q) + h(s).
\]
This shows that $\sigma_q=\sigma_p$ and the singularities on $\sigma=\sigma_p$ of $F(s,1/p)$ and $F(s,1/q)$ coincide, apart from a finite number of them. If $F\in\S$, in view of Lemma 7 this is possible only when the differences of the involved arithmetic progressions are the same, i.e. when $p=q$, a contradiction. If $F\in\S^\sharp$ we argue as in the proof of Theorem 1, using the idea in Gierszewski \cite{Gie/2020}. Thus $F_p(s)^{-1}$ has no singularities and hence it is a polynomial in $p^{-s}$. By symmetry of arguments the same is true for $F_q(s)^{-1}$, and the result follows. \qed

\section{Proof of Theorem 3}

\smallskip
Clearly, $F$ splits at every prime since it belongs to $\S$. Moreover, denoting by $p$ the prime number $q_F$, given a prime $p_1\neq p$ obviously there exists a distinct prime $p_2\neq p$ such that $p_1\equiv p_2$ (mod $p$). Hence, thanks to Theorem 2, we only have to show that $F_p(s)^{-1}$ is a polynomial in $p^{-s}$.

\smallskip
By orthogonality, for $\si>1$ we have
\begin{equation}
\label{eq:t3-1}
\begin{split}
\sum_{a=1}^{p-1} F(s,a/p) &= (p-1) \sum_{p|n} \frac{a(n)}{n^s} - \sum_{p\nmid n} \frac{a(n)}{n^s} \\
&= (p-1) \big(F(s) - F_p(s)^{-1}F(s)\big) - F_p(s)^{-1}F(s) \\
&= \left(p -1 -pF_p^{-1}(s)\right) F(s).
\end{split}
\end{equation}
For every $1\leq a<p$ we fix a prime  $p_a\equiv a$ (mod $p$), so that $F(s,a/p) =F(s,p_a/p)$. Thus by Lemma 4 with $\alpha = p_a\slash p$ and $K$ arbitrarily large we obtain
\begin{equation}
\label{eq:t3-2}
 F(s,a/p) = -i\omega_F^* \left(\frac{p_a}{\sqrt{p}}\right)^{2s-1+i\theta_F} \sum_{\nu=0}^K  \left(\frac{i p_a}{2\pi}\right)^{\nu} Q_{\nu}(s) \overline{F}(s+\nu+2i\theta_F, -1/p_a) +H_K(s,p_a/p).
 \end{equation}
But thanks to Lemma 1 we have
 \begin{equation}
 \label{eq:t3-3}
 \overline{F}(s, -1/p_a) = \frac{1}{1-p_a} \overline{F}_{p_a}(s)^{-1} \overline{F}(s)
 + \sum_{\substack{\chi \, {\rm (mod} \, p_a {\rm )}\\ \chi\neq\chi_0}}
 \overline{c(\chi,p_a)}\overline{F^\chi}(s).
 \end{equation}
From the hypothesis of Theorem 3 we have that the twists $\overline{F^\chi}(s)$ belong to $S^\sharp$ and have degree 2; moreover, by Theorem 2, $\overline{F}_{p_a}(s)^{-1}$ is a polynomial in $p_a^{-s}$. Therefore \eqref{eq:t3-2} and \eqref{eq:t3-3} give meromorphic continuation of $F(s,a\slash p)$ to the whole complex plane, with possible poles only at the points $1-\nu -i\theta_F$, $\nu\geq 0$. 
 
 \smallskip
In addition we have
\begin{equation}
\label{6-new1}
\overline{F}_{p_a}(s)^{-1}\ll p_a^{C|\sigma|}
\end{equation}
for a certain $C>0$ and, according to Lemma 3,
\begin{equation}
\label{6-new2}
\overline{F^\chi}(s)\ll A^{|\sigma|} |\sigma|^{2|\sigma|+1}
\end{equation}
for a certain $A>0$, uniformly for $\sigma\leq -1$ and $b\leq t\leq c$, with arbitrary fixed $b,c$ such that $|\theta_F|<b<c$. Thus, applying Lemma 5 with
\[ 
{\Q}_{\nu}(s) = \left(\frac{p_a}{2\pi}\right)^{\nu} |Q_{\nu}(s)| \qquad \text{and} \qquad {\F}(s) = |\overline{F}(s+2i\theta_F,-1/p_a)|, 
\]
we obtain that
\begin{equation}
\label{eq:t3-3a}
F(s,a/p) \ll B^{|\sigma|} |\sigma|^{2|\sigma|+1}
\end{equation}
for a certain $B>0$, uniformly for $\sigma\leq -1$ and $b\leq t\leq c$, with arbitrary fixed $|\theta_F| <b< c$. Now, rewriting 
\eqref{eq:t3-1} as
\begin{equation}
\label{eq:t3-4}
F_p(s)^{-1} = 1 - \frac{1}{p} - \frac{1}{pF(s)} \sum_{a=1}^{p-1} F( s,a/p),
\end{equation}
we obtain the meromorphic continuation of $F_p(s)^{-1}$ to the whole complex plane. Finally, using the periodicity and density arguments as in the proof of Theorem 1 we conclude that $F_p(s)^{-1}$  is an entire function. Moreover, by \eqref{eq:t3-3a}, \eqref{eq:t3-4} and Lemma 3 we have that
\[  
F_p(s)^{-1}\ll D^{|\sigma|}
\]
for certain $D>0$. Thus, again as in the proof of Theorem 1, $F_p(s)^{-1}$ is a polynomial in $p^{-s}$ and the result follows. \qed

\section{Proof of Theorem 4}

\smallskip
Similarly as in the proof of Theorem 3, for every $(a,q_F)=1$ we fix a prime $p_a\equiv a$ (mod $q_F$). From \eqref{t2-1} with $p=p_a$ we obtain
\[
F(s,1/p_a) \ll |F(s)| + |F_{p_a}(s)^{-1} F(s)| + \sum_{\substack{\chi \, {\rm (mod} \, p {\rm)} \\ \chi\neq \chi_0}} |c(\chi,p_a) F^{\chi}(s)|.
\]
As in the proof of Theorem 3 we have that \eqref{6-new1} and \eqref{6-new2} hold in the present case as well, the second bound with $F^\chi$ in place of $\overline{F^\chi}$ and $b>  \max_{\chi \, {\rm (mod} \, p {\rm)}} |\theta_{F^{\chi}}|$. Thus, thanks to Lemma 3, for such $s$ we have
\[
F(s,1/p_a) \ll C_1^{|\sigma|} |\sigma|^{2|\sigma|+1}
\]
with a certain $C_1>0$. Therefore, applying Lemma 4 with $\alpha=p_a/q_F$ and then Lemma 5, we obtain
\begin{equation}
\label{eq:t4-1}
F(s,p_a/q_F) \ll C_2^{|\sigma|}\sum_{0\leq\nu\leq|\sigma|+2} |Q_{\nu}(s)| |F(\overline{s}+\nu-2i\theta, 1/p_a)|
+ (A|\sigma|)^{|\sigma|} \ll C_3^{|\sigma|}|\sigma|^{2|\sigma|+1}
 \end{equation}
with certain $C_2,C_3>0$, uniformly for $\sigma\leq -1$ and $b<t<c$. 

\smallskip
 Let now $p$ be a sufficiently large prime, say $p\geq B_F$ with a certain $B_F>q_F$. Then there exists $(a,q_F)=1$ such that $p\equiv p_a$ (mod $q_F$). We shall estimate $|F_p(s)^{-1}|$ from above in the horizontal half-strip $\sigma\leq -1$, $b<t<c$, where $c> b+2\pi\slash \log p$ and $b$ is sufficiently large. Using Lemma 4 with $\alpha=1/p$, the fact that $F(s,p/q_F)= F(s,p_a/q_F)$, Lemma 5 and \eqref{eq:t4-1} we obtain
 \begin{equation}
 \label{eq:t4-2}
 \begin{split}
 F(s,1/p) &\ll   \left(\frac{p}{\sqrt{q_F}}\right)^{2|\sigma|+1} \sum_{0\leq\nu\leq |\sigma|+2} 
 \left(\frac{q_F}{2\pi p}\right)^{\nu} |Q_{\nu}(s) F(\overline{s}+\nu-2i\theta_F,p_a/q_F)|
 +(A|\sigma|)^{|\sigma|+1} \\
 &\ll (C_4 p)^{2|\sigma|+1} |\sigma|^{2|\sigma|+1}.
 \end{split}
\end{equation}
Next, from Lemma 2 with $m=1$ and then Lemma 3 and \eqref{eq:t4-2} we get, since the conductor $q_{F^\chi}$ of $F^\chi$ equals $q_Fp^2$ thanks to the hypotheses of Theorem 4, that
\[
\begin{split}
F_p(s)^{-1} &\ll \frac{1}{|F(s)|}\sum_{\substack{\chi \, {\rm (mod} \, p {\rm)} \\ \chi\neq \chi_0}} |c(\chi,p) F^{\chi}(s)| +
\frac{1}{|F(s)|} |F(s,1/p)| +1 \\
&\ll (C_5 p)^{2|\sigma|+1}\ll p^{(5\slash 2)|\sigma|}.
\end{split}
\]
But we already know from Theorem 2 that $F_p(s)^{-1}=P_p(p^{-s})$ for certain polynomial $P_p\in {\mathbb C}[z]$, hence the last inequality shows that its degree is at most $2$. Since the constant term of $P_p$ is $1$,  we can write
\[F_p(s) = \left(1-\frac{\alpha_p}{p^s}\right)^{-1}\left(1-\frac{\beta_p}{p^s}\right)^{-1}\]
for certain complex numbers $\alpha_p$ and $\beta_p$. Recalling that the  Dirichlet coefficients of $F$ satisfy Ramanujan's conjecture we must have $|\alpha_p|, |\beta_p| \leq 1$ (see p.448-449 of \cite{Ka-Pe/2015}), and the result follows. \qed

\bigskip
\bigskip

\ifx\undefined\bysame{poly}.
\newcommand{\bysame}{\leavevmode\hbox to3em{\hrulefill}\ ,}
\fi

\bigskip
\bigskip
\noindent
Jerzy Kaczorowski, Faculty of Mathematics and Computer Science, A.Mickiewicz University, 61-614 Pozna\'n, Poland and Institute of Mathematics of the Polish Academy of Sciences, 00-956 Warsaw, Poland. e-mail: \url{kjerzy@amu.edu.pl}

\medskip
\noindent
Alberto Perelli, Dipartimento di Matematica, Universit\`a di Genova, via Dodecaneso 35, 16146 Genova, Italy. e-mail: \url{perelli@dima.unige.it}

\end{document}